\documentclass[reqno]{amsart}
\usepackage{enumerate}
\usepackage{bbm}
\usepackage[margin=1.25in]{geometry}
\usepackage{amsmath,amsthm,amscd,amssymb}
\usepackage{latexsym}
\usepackage{color}
\usepackage{todonotes}

\numberwithin{equation}{section}

\newcommand{\N}{{\mathbb N}}
\newcommand{\XX}{{\mathfrak X}}
\newcommand{\X}{{\bf X}}
\newcommand{\T}{{\mathbf T}}
\newcommand{\1}[1]{{\mathbbm{1}\left( #1\right)}}
\renewcommand{\P}{{\mathbf P}}

\newcommand{\p}{{\mathcal P}}
\newcommand{\E}{{\mathcal E}}

\theoremstyle{plain}
\newtheorem{theorem}{Theorem}
\newtheorem{corr}[theorem]{Corollary}
\newtheorem{lemma}[theorem]{Lemma}
\newtheorem{prop}[theorem]{Proposition}
\theoremstyle{definition}
\newtheorem{defi}[theorem]{Definition}

\theoremstyle{remark}

\title[Markovian stick-breaking priors]{On the use of Markovian stick-breaking priors}

\author{William Lippitt}
\address{Mathematics, University of Arizona, Tucson, AZ 85721}
\email{wlippitt@math.arizona.edu}

\author{Sunder Sethuraman}
\address{Mathematics, University of Arizona, Tucson, AZ 85721}
\email{sethuram@math.arizona.edu}

\begin{document}

\begin{abstract}
In \cite{DLS}, a `Markovian stick-breaking' process which generalizes the Dirichlet process $(\mu, \theta)$ with respect to a discrete base space $\XX$ was introduced.  In particular, a sample from from the `Markovian stick-breaking' processs may be represented in stick-breaking form $\sum_{i\geq 1} P_i \delta_{T_i}$ where $\{T_i\}$ is a stationary, irreducible Markov chain on $\XX$ with stationary distribution $\mu$, instead of i.i.d. $\{T_i\}$ each distributed as $\mu$ as in the Dirichlet case, and $\{P_i\}$ is a GEM$(\theta)$ residual allocation sequence.  Although the motivation in \cite{DLS} was to relate these Markovian stick-breaking processes to empirical distributional limits of types of simulated annealing chains, these processes may also be thought of as a class of priors in statistical problems.  The aim of this work in this context is to identify the posterior distribution and to explore the role of the Markovian structure of $\{T_i\}$ in some inference test cases.
\end{abstract}



\maketitle

\begin{center}
{\large \em Dedicated to Professor M.M. Rao on his 90th birthday.}
\end{center}

\section{Introduction}

Let $\XX\subseteq\N$ be a discrete space, either finite or countable.  Let also $\mu$ be a measure on $\XX$, and $\theta>0$ be a parameter.  The Dirichlet process on $\XX$ with respect to pair $(\theta, \mu)$ is an object with fundamental applications to Bayesian nonparametric statistics (cf. books \cite{Ghosal_VanVaart}, \cite{Muller}).  Formally, the Dirichlet process is a probability measure on the space of probability measures on $\XX$ such that a sample $\mathcal{P}$, with respect to any finite partition $(A_1,\ldots, A_k)$ of $\XX$, has the property that the distribution of $(\mathcal{P}(A_1),\ldots, \mathcal{P}(A_k))$ is Dirichlet with parameters $(\theta \mu(A_1),\ldots, \theta \mu(A_k))$ (cf. \cite{Ferguson}, \cite{Blackwell_McQueen}).  Importantly, the Dirichlet process has a `stick-breaking' representation:  A sample $\mathcal{P}$ can be represented in form $\sum_{j=1}^\infty P_j \delta_{T_j}$ where $\P = \{P_j\}_{j\geq 1}$ is a GEM$(\theta)$ residual allocation sequence, and $\{T_j\}_{j\geq 1}$ is an independent sequence of independent and identically distributed (i.i.d.) random variables on $\XX$ with common distribution $\mu$ (cf. \cite{JS}, \cite{Pitman}).  Here, a GEM sequence is one where $P_1 = X_1$ and $P_j = X_j \big(1-\sum_{i=1}^{j-1}P_i\big)$ for $j\geq 2$, and $\{X_j\}_{j\geq 1}$ are i.i.d. Beta$(1,\theta)$ random variables.

There are several types of generalizations of the Dirichlet process in the literature such as Polya tree, and species sampling processes \cite{Ghosal_VanVaart}[Ch. 14], \cite{Lavine}, among others.  In \cite{DLS}, another generalization where $\{T_j\}_{j\geq 1}$ is a Markov chain was introduced:  Let $G=\{G_{i,j}: i,j\in \XX\}$ be a generator matrix, that is $G_{i,j}\geq 0$ for $i\neq j$ and $G_{i,i} = -\sum_{j\neq i} G_{i,j}$, which is irreducible with suitably bounded entries, and has $\mu$ as its stationary distribution.  Let also $Q= I + G/\theta$ be a Markov transition kernel on $\XX$ with stationary distribution $\mu$.  Now, define $\T=\{T_j\}_{j\geq 1}$ as the stationary Markov chain with transition kernel $Q$.  The `Markovian stick-breaking' process is then represented as $\sum_{j\geq 1}P_j \delta_{T_j}$ where again $\P$ is an independent GEM$(\theta)$ sequence.

Although it was shown in \cite{DLS}, \cite{DS} that such Markovian stick-breaking processes connect to the limiting empirical distribution of certain simulated annealing chains, it is natural to consider their use as priors in statistical problems, the aim of this article.  We first give a formula for the moments of the Markovian stick-breaking process in Theorem \ref{thmpriornset}.  Then, we compute the posterior distribution moments in terms of this formula in Proposition \ref{propposterior}, Corollary \ref{cor_post}.  Consistency of the posterior distribution is stated in Proposition \ref{prop_consistency}, noting the full support property of the process in Proposition \ref{prop_support}.  In Proposition \ref{prop_theta}, we discuss asymptotics of the process with respect a `strength' parameter.

A main part of this work is also to consider the use and behavior of the Markovian stick-breaking process in as a prior for inference of histograms. In this context, the generator $G$ can be thought as a priori belief of weights or affinities in a `network' of categories.  For instance, in categorical data, one may believe that affinities between categories differ depending on the pair, and also that they may be directed in hierarchical situations.  In using Dirichlet priors, there is an implicit assumption that the network connecting categories is complete and the affinity between two categories cannot depend on both categories.  However, in using a Markovian stick-breaking prior, one can build into the prior a belief about the weight structure on the network by specification of the generator $G$.  In simple experiments, we show interesting behaviors of the posterior distribution from these Markovian stick-breaking priors, in comparison to Dirichlet priors.

The structure of the paper is to define carefully the Markovian stick-breaking process in Section \ref{def_section}.  Then, in Section \ref{results_sect}, we state results on their moments, posterior distribution, and consistency.  In Section \ref{use_section}, we discuss the use of these processes as priors, present simple numerical experiments, and provide some context with previous literature.

\section{Definition of the Markovian stick-breaking process}
\label{def_section}

	We take as convention empty sums are 0, empty products of scalars are 1, empty products of matrices are the identity, and that a product of matrices is computed as
	$$\prod_{j=1}^n M_j=M_n\cdot M_{n-1}\cdot\ \cdots\ \cdot M_1.$$
	For a set $A\subseteq\XX$, we define $D(A)$ as the diagonal square matrix over $\XX$ with entries $D_{xx}(A)=\delta_x(A)$.   For $x\in \XX$, let $e_x\in {\mathbb R}^\XX$ be the column vector with a $1$ in the $x$th entry and $0$'s in other entries.  Let also $\vec 1$ be the column vector of all $1$'s.

We now give a precise definition of the well-known GEM `Griffiths-Engel-McCloskey' residual allocation sequence, which apportions a unit resource into infinitely many parts.
\begin{defi}[GEM]
	Let $\X=(X_j)_{j=1}^\infty$ be an i.i.d. collection of Beta$(1,\theta)$ variables for some positive constant $\theta$.  Define $\P=(P_j)_{j=1}^\infty$ by $P_1=X_1$ and $P_j = X_j\big(1-\sum_{i=1}^{j-1}P_i\big)$ for $j\geq 2$, which leads to the formula
	$$P_j=X_j\prod_{i=1}^{j-1}(1-X_i).$$
	We say $\P$ has GEM$(\theta)$ distribution.
\end{defi}

To define the `Markovian stick-breaking' process on the discrete space $\XX$, we now state carefully the definition of a Generator kernel or matrix.
\begin{defi}[Generator]
We call a real-valued matrix $G=(G_{xy})_{x,y\in\XX}$, a generator matrix over $\XX$, if
	\begin{enumerate}
		\item For each pair $x,y\in\XX$ with $x\neq y$, then $G_{x,y}\ge0$.
		\item For each $x\in\XX$, $G_{xx}=-\sum_{y\in\XX-\{x\}}G_{xy}$.
		\item $\theta^G:=\sup_{x\in\XX}|G_{xx}|<\infty$
	\end{enumerate}
	If $\mu$ is a stochastic vector over $\XX$ and $\mu^TG=\vec 0$, we call $\mu$ a stationary distribution of $G$.
\end{defi}

	Note that if $\theta\ge\theta^G$ for a generator matrix $G$, then 
	$$Q=I+G/\theta,$$
	where $I$ is the identity kernel, is a stochastic matrix over $\XX$. All such $Q$'s share stationary distributions and communication classes. As such, we refer to the stationary distributions and irreducibility properties of $G$ and $Q$'s interchangeably.

We now define the `Markovian stick-breaking' measure (MSB) as follows.
\begin{defi}[MSB$(G)$]\label{defnMSBM}
	Let $G$ be an irreducible positive recurrent generator matrix over a discrete space $\XX$, with stationary distribution $\mu$.  Let $\theta\ge\theta^G$, and stochastic matrix $Q=I+G/\theta$.  Let $\P\sim$ GEM$(\theta)$, and let $\T=(T_j)_{j=1}^\infty$ be a stationary, homogeneous Markov chain in $\XX$ with transition kernel $Q$ and independent of $\P$. Define the random measure $\nu$ over $\XX$ by
	$$\nu=\sum_{j=1}^\infty P_j\delta_{T_j}\ .$$
	We say $\nu$ has Markovian stick-breaking distribution with generator $G$, and the pair $(\nu,T_1)$ has MSB$(G)$ distribution.  \end{defi}

We note, in this definition, the distribution of $\nu$ does not depend on the choice of $\theta\geq \theta^G$, say as its moments by Corollary \ref{mom_corr} below depend only on $G$; see also \cite{DLS} for more discussion.
We remark also, when $Q$ is a `constant' stochastic matrix with common rows $\mu$, then $\{T_j\}_{j\geq 1}$ is an i.i.d. sequence with common distribution $\mu$ and so the Markovian stick-breaking measure $\nu$ reduces to the Dirichlet distribution with parameters $(\theta, \mu)$; see \cite{DLS} for further remarks.

\section{Results on moments, posterior distribution, and consistency}
\label{results_sect}

We now compute in the next formulas certain moments of the Markovian stick-breaking measure with respect to generator $G$, which identify the distribution of $\nu$.

\begin{theorem}
	Let $G$ be an irreducible positive recurrent generator matrix on $\XX$, and let $(\nu,T_1)\sim$ MSB$(G)$.  Let also $(A_j)_{j=1}^n$ be a collection of disjoint subsets of $\XX$, $x\in\XX$, and $\vec k\in\{0,1,2,\ldots\}^n$.  Then,
	$$\E\left[{\prod_{j=1}^n\nu(A_j)^{k_j}\Big|T_1=x}\right]=\left(\#\mathbb{S}(\vec k)\right)^{-1}\sum_{\sigma\in\mathbb{S}(\vec k)}e_x^T\left[\prod_{j=1}^k\left((I-G/j)^{-1}D(A_{\sigma_j})\right)\right]\vec 1$$
	where $k=\sum_{j=1}^n k_j$, $\mathbb S(\vec k)$ is the collection of distinct permutations of $k$-lists of $k_1$ many 1's, $k_2$ many 2's, and so on to $k_n$ many $n$'s, and $\#\mathbb S(\vec k)$ is the cardinality of this set.  \label{thmpriornset}
\end{theorem}

The proof of Theorem \ref{thmpriornset} is given in the Section \ref{proof_sec}.

\begin{corr} 
\label{mom_corr}
In the context of the previous theorem, suppose $A_j=\{x_j\}$.  Then,
	$$\E\left[{\prod_{j=1}^n \nu(x_j)^{k_j}}\right]=\left(\#\mathbb{S}(\vec k)\right)^{-1}\sum_{\sigma\in\mathbb{S}(\vec k)}\mu_{x_{\sigma_k}}\prod_{j=1}^{k-1}(I-G/j)^{-1}_{x_{\sigma_{j+1}},x_{\sigma_j}}$$
\end{corr}

We remark that Corollary \ref{mom_corr} is an improvement of a corresponding formula in \cite{DS} found, by different means, when $\XX$ is finite and $G$ has no nonzero entries.
\medskip

These formulas will be of help to identify the posterior distribution, if the Markovian stick-breaking measure is used as a prior.  In the case of the Dirichlet process, the posterior distribution is again in the class of Dirichlet processes:  Namely given $\nu$, let $Y_1,\ldots, Y_n$ be i.i.d. random variables with distribution $\nu$.  Then, the distribution of $\nu$ given $Y^n=\{Y_j\}_{j=1}^n$ is a Dirichlet process with parameters $(\theta, \mu + \sum_{j=1}^n \delta_{Y_j})$.  However, when $\nu$ is a general Markovian stick-breaking measure, such a neat correspondence is not clear.  But, later in Proposition \ref{propposterior}, we write the posterior moments in terms of `size-biased' moments with respect to the prior.

We now give a representation of a sequence $Y_1,\ldots, Y_n$, conditional on a sample $\nu$ from the Markovian stick-breaking process, which is i.i.d. with common distribution $\nu$.  This representation is standard with respect to the Dirichlet process and relatives such as species sampling processes (cf. Ch. 14 \cite{Ghosal_VanVaart}).

\begin{prop}\label{propdata}
	Consider the Markovian stick-breaking process $\nu$ built from $\P$ and $\T$.  For $n\geq 1$, let $(J_i)_{i=1}^n$ be a collection of positive integer valued random variables such that $\p \big(J_i=j_i:1\le i\le n\big|\P,\T\big)=\prod_{i=1}^n P_{j_i}$.
	Define the sequence $Y^n=(Y_i)_{i=1}^n$ where $Y_i=T_{J_i}$ for $1\leq i\leq n$.   Then, $Y^n\Big|\nu,T_1$ is a collection of i.i.d. variables taking values in $\XX$ with common distribution $\nu$.
\end{prop}

\begin{proof} Compute, noting $\nu(x)=\sum_{j\geq 1}P_j \1{T_j = x}$, that
\begin{align*}
		&\p\big(Y^n=y^n\big|\P,\T\big)\  = \ \sum_{j_1=1}^\infty\cdots\sum_{j_n=1}^\infty\p\big(J_i=j_i,T_{j_i}=y_j:1\le i\le n\big|\P,\T\big)\\
		&\ \  =\sum_{j_1=1}^\infty\cdots\sum_{j_n=1}^\infty\prod_{i=1}^nP_{j_i}\1{T_{j_i}=y_i}\ = \ 
		\prod_{i=1}^n\sum_{j=1}^\infty P_j\1{T_j=y_i}\ = \ 
		\prod_{i=1}^n\nu(y_i)
	\end{align*}
	Since $(\nu,T_1)$ is a function of $\P$ and $\T$ and $\p(Y^n = y^n|\nu, T_1) = \E[\p\big(Y^n=y^n\big|\P,\T\big)|\nu, T_1]$, the result follows.
\end{proof}

The following identifies the posterior distribution in terms of its moments, given as certain `size-biased' expressions with respect to the prior.

\begin{prop}\label{propposterior}
	Let $\nu$ be a random probability measure taking values in the simplex $\Delta_\XX=\{ {\bf p}: \sum_{x\in\XX} p_x = 1, 0\leq p_x\leq 1\}$, and let $T$ be an $\XX$-valued random variable. For $n\geq 1$, let $Y^n=(Y_1,...,Y_n)$ be a sequence of random variables such that $Y^n\Big|\nu,T$ are i.i.d. with common distribution $\nu$. 
	Let also $\vec k = \vec k^n$ denote the frequencies of $y^n$, that is, for each $x\in\XX$, $k_x=\#\{j:1\le j\le n,\ y_j=x\}$.  In addition, let $\vec l\in\{0,1,2,\ldots\}^\XX$, where $l_x =0$ except for finitely many $x\in \XX$.  
	Then, for events $A\in\sigma(T)$, we have
	$$\E\left[\prod_{x\in\XX}\nu(x)^{l_x}\Big|Y^n=y^n,A\right]=\frac{\E\Big[\prod_{x\in\XX}\nu(x)^{k_x+l_x}\Big|A\Big]}{\E\Big[\prod_{x\in\XX}\nu(x)^{k_x}\Big|A\Big]}$$
\end{prop}

\begin{proof}
Define $m$ additional random variables $Y_{n+1},Y_{n+2},...,Y_{n+m}$, by augmenting the probability space if necessary, such that together $Y_1,Y_2,\ldots, Y_{n}, Y_{n+1},\ldots, Y_{n+m}| \nu, T$ are i.i.d. with common distribution $\nu$. In particular, for $y^{n+m}\in \XX^{n+m}$, we have $\p\big(Y^{n+m}=y^{n+m}\big|\nu,T\big)=\prod_{j=1}^{n+m}\nu(y_j)$.
	
	Recall now $\vec k$ the frequencies of $y^n$. Let 
	$(y_{n+1},y_{n+2},...,y_{n+m})\in\XX^m$ be any sequence with frequencies $\vec l$. Then, we compute
	\begin{align*}
		&\E\left[\prod_{x\in\XX}\nu(x)^{l_x}\Bigg|Y^n=y^n,A\right] \  = \ \E\left[\frac{\p\big(Y^{n+m}=y^{n+m}\big|\nu,T\big)}{\p\big(Y^n=y^n\big|\nu,T\big)}\Big |Y^n=y^n,A\right]\\
		 & =\E\Big[\p\big(Y^{n+m}=y^{n+m}\big|Y^n=y^n,\nu,A\big)\Big|Y^n=y^n,A\Big]\ = \ 
		  \p\big(Y^{n+m}=y^{n+m}\big|Y^n=y^n,A\big)\\
		 & =\frac{\p\big(Y^{n+m}=y^{n+m}\big|A\big)}{\p\big(Y^n=y^n\big|A\big)}\ = \ 
		 \frac{\E\Big[\prod_{x\in\XX}\nu(x)^{k_x+l_x}\Big|A\Big]}{\E\Big[\prod_{x\in\XX}\nu(x)^{k_x}\Big|A\Big]},
	\end{align*}
	and the result follows. \end{proof}

Returning to the Markovian stick-breaking process $\nu$, 
given the `data' $Y_1,\ldots, Y_n$ conditional on $\nu$ and $T_1$, we may evaluate the posterior moments of the Markovian stick-breaking measure as a case of Proposition \ref{propposterior}.

\begin{corr} 
\label{cor_post}
Let $\nu$ be a Markovian stick-breaking process, and $Y_1,\ldots, Y_n| \nu, T_1$ be i.i.d. with common distribution $\nu$ (say, as in Proposition \ref{propdata}). Let also $\vec k=(k_x)_{x\in\XX}$ be such that $k_x=\#\{j:1\le j\le n\ and\ y_j=x\}$ for each $x\in \XX$. Let in addition $\vec l\in\{0,1,2,\ldots\}^\XX$ be a vector with only finitely many non-zero entries, and $m=\sum_{x\in\XX}l_x$. Then, for each $x\in\XX$, we have
	\begin{align*}
		\E\left[\prod_{w\in\XX}\nu(w)^{l_w}\Big|Y^n=y^n,T_1=x\right] 
		& =\frac{\left(\#\mathbb S(\vec k)\right)\hspace{-.2mm}\sum\limits_{\sigma\in\mathbb S(\vec k+\vec l)}\hspace{-.2mm}(I-G/(n+m))^{-1}_{x,\sigma_{n+m}}\prod\limits_{j=1}^{n+m-1}(I-G/j)^{-1}_{\sigma_{j+1},\sigma_j}}
		{\left(\#\mathbb S(\vec k+\vec l)\right)\sum\limits_{\sigma\in\mathbb S(\vec k)}(I-G/n)^{-1}_{x,\sigma_n}\prod\limits_{j=1}^{n-1}(I-G/j)^{-1}_{\sigma_{j+1},\sigma_j}}\\
		\\
	{\rm and \ \ }	\E\left[\prod_{w\in\XX}\nu(w)^{l_w}\Big|Y^n=y^n\right] & =\frac{\left(\#\mathbb S(\vec k)\right)\sum\limits_{\sigma\in\mathbb S(\vec k+\vec l)}\mu_{\sigma_{n+m}}\prod\limits_{j=1}^{n+m-1}(I-G/j)^{-1}_{\sigma_{j+1},\sigma_j}}{\left(\#\mathbb S(\vec k+\vec l)\right)\sum\limits_{\sigma\in\mathbb S(\vec k)}\mu_{\sigma_n}\prod\limits_{j=1}^{n-1}(I-G/j)^{-1}_{\sigma_{j+1},\sigma_j}}.
	\end{align*}
\end{corr}

We now give a statement of `consistency' with respect to the posterior distribution, in line with limits of `Bayes estimators' in \cite{Freedman}, by considering the moment expression in Proposition \ref{propposterior}, when $\XX$ is finite.   
Consistency, in the case $\XX$ is countably infinite, may be pathological according to \cite{Freedman}, and so we limit out discussion accordingly.

\begin{prop}
\label{prop_consistency}
Let $\nu$ be a Markovian stick-breaking process on a finite state space $\XX$.  Let also $Y_1,\ldots, Y_n|\nu, T_1$ be i.i.d. with common distribution $\nu$.  

Suppose, for each $x\in \XX$ as $n\uparrow\infty$, that
$\frac{1}{n}\sum_{j=1}^n \1{Y_j = x} \rightarrow \eta_x$ a.s.,
where $\eta = \{\eta_x\}_{x\in \XX}\in \Delta_\XX$. Then, as $n\uparrow\infty$, the posterior distribution 
$\mu_n = \P(\nu\in \cdot |Y^n)$
converges a.s. to $\delta_\eta$.
\end{prop}

\begin{proof}
Write,
$\mu_n(B) = \P(\nu\in B| Y^n=y^n) = \frac{\P(Y^n=y^n, \nu\in B)}{\P(Y^n=y^n)} = \frac{\E\big[\prod_{j=1}^n \nu(y_j), \nu\in B\big]}{\E\big[\prod_{j=1}^n \nu(y_j)\big]}$.
By Theorem 1 in \cite{Freedman}, if $\eta$ belongs to the support of $\nu$, the desired convergence of $\mu_n$ to $\delta_\eta$ follows.

Hence, to finish, we note by Proposition \ref{prop_support} below that $\nu$ has full support on the simplex $\Delta_\XX$. 
\end{proof}

The following is an improvement of a corresponding result in \cite{DS} when $G$ has no zero entries, by directly considering the stick-breaking form of $\nu$.

\begin{prop}
\label{prop_support}
For finite $\XX$, the Markovian stick-breaking measure $\nu$ with respect to irreducible $G$ has full support on the simplex $\Delta_\XX$.
\end{prop}

\begin{proof}

Let $r=|\XX|$. Since $\nu$ has the form $\sum_{j\ge 1}P_j\delta_{T_j}$, the idea is to consider a path of the Markov chain $\bf T$ with prescribed visits to states 1,2,...,$r$, and realizations of the GEM$(\theta)$ sequence $\bf P$ with values such that $\nu$ belongs to a small $\epsilon$-ball around $\eta$. 

Since $Q$ is irreducible, there exists an integer $n\ge r$ and a path $(t_1,...,t_n)\in\XX^n$ such that the chain $\bf T$ has positive probability of starting on the path $P(T_i=t_i:1\le i\le n)>0$ and such that the path hits every state $x\in\XX$. For each state $x\in\XX$, define $i_x=\min\{i:t_i=x\}\in\{1,...,n\}$ to be the first time the path hits state $x$. 

Since $\bf P$ is distributed as a residual allocation model constructed from iid proportions each having full support on the unit interval, $\bf P$ has full support on $\Delta_\infty$. Thus, for each $\delta>0$, we have with positive probability that simultaneously $P_{i_x}>\eta_x-\delta$ for all $x\in\XX$. Noting the following containment of events 
$$\big\{1\le i\le n:T_i=t_i;\ \forall x: P_{i_x}>\eta_x-\delta\big\}\subseteq \big\{\forall x:-\delta<\nu_x-\eta_x<(r-1)\delta\big\}\subseteq \big\{\forall x:|\nu_x-\eta_x|<(r-1)\delta\big\}$$
and taking $\delta=\epsilon/(r-1)$, we then have
$$\p(\forall x\in\XX:|\nu_x-\eta_x|<\epsilon)\ge \p(1\le i\le n:T_i=t_i;\ \forall x\in\XX:P_{i_x}>\eta_x-\epsilon/(r-1))>0.$$
Hence, $\nu$ is within $\epsilon$ of $\eta$ with positive probability.
\end{proof}

When the stochastic matrix $Q$ is fixed, the parameter $\theta$ in the representation of $G = \theta(Q-I)$ can be viewed as a type of `strength' of the Markovian stick-breaking $\nu$, as more discussed in the next section.

\begin{prop}
\label{prop_theta}
Let $Q$ be irreducible positive recurrent stochastic and define $G^\theta=\theta(Q-I)$.  Then, the Markovian stick-breaking measure $\nu=\nu^{(\theta)}$ parametrized by $G^\theta$ converges in probability to the stationary vector $\mu$ of $Q$ as $\theta\uparrow\infty$.
\end{prop}

\begin{proof}Suppose $Q$ is aperiodic. Let $x\in \XX$. For all $\theta>0$ we have $E[\nu(x)]=\mu_x$.  As $\theta\uparrow\infty$, we have by Cor. 4.1 for each $n\in\N$ that
\begin{align*}
\lim_{\theta\rightarrow\infty}\E[\nu(x)^2] & =(I-G^\theta)^{-1}_{xx}\mu_x=\lim_{\theta\rightarrow\infty}\frac{\mu_x}{\theta+1}\sum_{j=0}^{n-1}\left(\frac{\theta Q}{\theta+1}\right)^j_{xx}+\frac{\mu_x}{\theta+1}\left[\left(\frac{\theta Q}{\theta+1}\right)^n\sum_{j=0}^{\infty}\left(\frac{\theta Q}{\theta+1}\right)^j\right]_{xx}\\
 & =\lim_{\theta\rightarrow\infty}0+\mu_x\left[Q^n(I-G^\theta)^{-1}\right]_{xx}=\mu_x^2
\end{align*}
since $Q^n$ converges to a constant stochastic matrix with rows $\mu$ as $n\rightarrow\infty$ and $(I-G^\theta)^{-1}$ is stochastic. Therefore, $\nu(x)$ converges in probability to $\mu_x$ as $\theta\uparrow\infty$.

If $Q$ is periodic, define aperiodic $Q'=0.5(Q+I)$ and note $G^\theta=\theta(Q-I)=2\theta(Q'-I)$. Since the proposition has been shown to apply to $Q'$, the result holds also for $Q$. \end{proof}

\section{On use of the MSB$(G)$ measure as a prior}
\label{use_section}

We explore in this section the use of the Markovian stick breaking measure MSB$(G)$ as a prior for multinomial probabilities.  In a nutshell, with respect to such a prior, when $G$ is in form $G=\theta(Q-I)$, the matrix $Q$ specifies an affinity network which reflects prior beliefs of association among categories.  Given observed data, the posterior mean histogram then computed will have the effect of `smoothing' the empirical probability mass function (pmf) according to the affinity network, in that mass levels of related categories will tend be similar.  The parameter $\theta$ as we will note will then represent a relative strength of this `smoothing'. 
 In particular, we consider, in simple examples, effects on the posterior mean histograms with respect to a few MSB$(G)$ priors in relation to Dirichlet priors, which do not assert affinities among categories.

Of course, `histogram smoothing' in the context of pmf estimation is an old subject with several Bayesian approaches.  For instance, see Leonard \cite{Leonard}, where multivariate logistic-normal priors are considered; Dickey and Jiang \cite{DJ}, where `filtered' Dirichlet distributions are proposed; Wong \cite{Wong}, where generalized Dirichlet distributions are used; and more recently Demirhan and Demirhan \cite{Demirhan}; see also the survey Agresti and Hitchcock \cite{AH}, and books Agresti \cite{Agresti}, Ghosal and Van der Vaart \cite{Ghosal_VanVaart}, and Congdon \cite{CongdonCat,CongdonMod} and references therein.
We remark there is also a large body of work for `histogram smoothing' with respect to Bayesian density estimation for continuous data, not unrelated to that for pmf inference.  See, for instance, Petrone \cite{Petrone}, Escobar and West \cite{Escobar}, and Hellmayr and Gelfand \cite{HG}, and references therein.

Similarly, categorical data may be viewed in terms of contingency tables with prior beliefs that certain factor outcomes are likely to co-occur or to occur separately, or that outcomes are likely to share a majority of factors.  Again, there is considerable work on Bayesian inference in this vein.  For instance, see Agresti and Hitchcock \cite{AH}, and books Agresti \cite{Agresti}, Ghosal and Van der Vaart \cite{Ghosal_VanVaart}, and Congdon \cite{CongdonCat,CongdonMod} and references therein.

\medskip\noindent \textbf{Histogram smoothing: Toy problem.}  We recall informally a basic `toy problem', with respect to the inference of the distribution of say shoe sizes, to set-up the main ideas.  Suppose a shoe seller is opening a new shop in town and wants to know the distribution of shoe sizes of the town population before stocking the shelves. Suppose that a person's shoe size is determined by their foot length, and that foot lengths are approximately Normal in distribution. Then, of course, we would expect that a histogram of shoe sizes would look approximately like a binned Normal histogram.

The shoe seller records the shoe sizes from a sample of individuals in town. In this multinomial data, categories are shoe sizes. We have some prior understanding of the context. Shoe sizes have a lower and upper bound, and presumably most people have shoe sizes relatively in the middle. Moreover, prior knowledge that shoe sizes arise from a continuous Normally distributed factor (foot length), would indicate that gaps in the shoe size sample histogram are likely not present in the true histogram.

One could use a Dirichlet prior, conveniently conjugate with multinomial data, though we will see shortly that such a prior cannot encompass all of the prior knowledge. Suppose there are $d$ possible shoe sizes/categories, numbered $1,2, \ldots, d$.  We specify a Dirichlet prior with parameters $(\theta \mu)$ where $\mu\in\Delta_d$ is the best guess at the shoe size probability mass function, and concentration parameter $\theta>0$ represents the level of confidence in the best guess $\mu$. If there is no `best guess', one could take $\mu = (1/d)(1,\ldots, 1)$, the uniform stochastic vector, and $\theta$ small. 
 
Let $\vec f\in\{0,1,2,\ldots\}^d = (f_1, f_2,\ldots, f_d)$ be the count vector from the sample of size $n$ collected, where $f_i$ is the number of people in the sample with shoe size $i$.  With this data in hand, one updates the prior belief by computing the posterior distribution.  In the case, if the prior is Dirichlet$(\theta\mu)$, the posterior would be Dirichlet$(\theta\mu + \vec f)$.  Then, the posterior estimate of the population distribution of shoe sizes would be the posterior mean $(\theta\mu + \vec f)/(\theta + n)$.

As an example, consider sample shoe size data collected from $15$ Normal samples binned into $d=16$ shoe sizes. Suppose we specify a so-called non-informative prior with $\mu=(1/d)\vec 1$ and $\theta=4$. In Figure 1, we see the prior estimate of the pmf in the left plot (i.e. $\mu$) represented as a histogram. In the middle plot is the empirical pmf computed from the sample. The posterior mean histogram, a weighted average of the left and middle plots, is seen in the right plot.

\begin{figure}
\label{Fig1}
\begin{center}
\includegraphics
{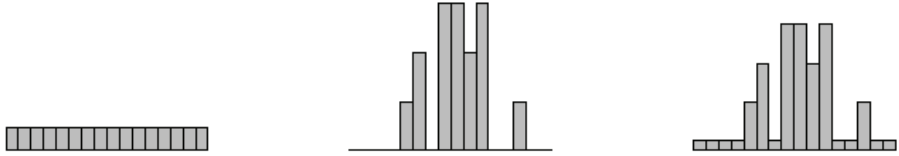}
\end{center}
\caption{Left plot is the prior mean histogram for a Dirichlet prior with uniform mean; Middle plot is of data collected; Right plot is the posterior histogram.}
\end{figure}

The posterior mean histogram is `smoother,' or less jagged, in that the two gaps in the data histogram have been partially filled. However, a Dirichlet prior does not allow too much control:  There is no notion of association between categories built in to the prior.  As such, one wouldn't be able to impose in some way that an empty bin between two `tall' bars should be filled with a similarly `tall' bar, or that an empty bin very far from any observed data should be left approximately empty.

In this context, we explore now use of a Markovian stick-breaking prior, which encodes associations between categories through specification of a network represented by the generator matrix $G$.  In this general network, categories are nodes and edges, directed or undirected, specify affinity between categories. The adjacency matrix for this network is then formed into the generator matrix $G$ by modifying diagonal entries appropriately to create generator matrix structure.  Recall that the matrix $G$, in the form $G=\theta(Q-I)$, specifies the transition matrix $Q$ for the Markovian sequence $\T$ as well as the parameter $\theta$ for the GEM sequence $\P$.  Accordingly, counts in the different categories are associated not only with respect to the GEM $\P$ but also with respect to the Markovian $\T$. 

We recall, in the Dirichlet$(\theta\mu)$ context, where $\T$ is an i.i.d. sequence with common distribution $\mu$ and $\P$ is GEM$(\theta)$, that the parameter $\theta$ is viewed as a `strength', and can represent in a sense the number of data points equivalent to the prior `belief'.  The corresponding posterior mean mass function is the weighted average $(\theta \mu + \vec{f})/(\theta + n)$ where $(1/n)\vec{f}$ is the empirical data probability mass function.  When $\theta\uparrow\infty$, the limit is the prior belief mean $\mu$.  
 
It is similar in the Markovian stick-breaking setting: If say the transition matrix $Q$ representing the network is specified in advance, the parameter $\theta$ is also a sort of relative strength in that, as $\theta\uparrow\infty$, $\nu$ converges in probability to $\mu$ (Proposition \ref{prop_theta}).

\begin{figure}
\label{Fig2}
\begin{center}
\includegraphics
{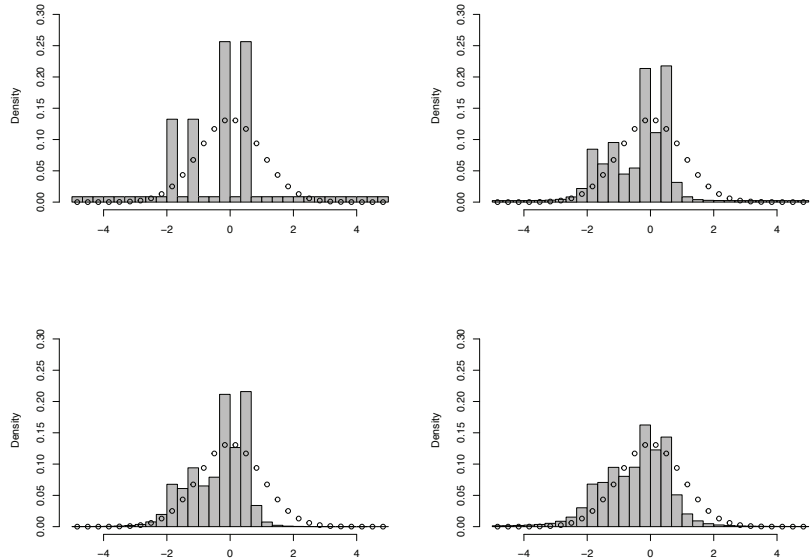}
\end{center}
\caption{Normal$(0,1)$ population pmf given by dotted curve; $6$ samples across $30$ bins in range $[-5,5]$ (not pictured) with $1$ point in bins $10, 12$ and $2$ points in bins $15, 17$.  Posterior mean mass functions from MSB$(G)$ priors:  Top left: $G_1$=Dirichlet$(w,\ldots, w)$, $w=2/29$; Bottom left: $G_2=$Tridiagonal with $w=3$; Bottom right: $G_3=$Tridiagonal with $w=8$; Top right $G_4= (G_1 + 2.5G_2)/3.5$.}
\end{figure}

\begin{figure}
\label{Fig3}
\begin{center}
\includegraphics
{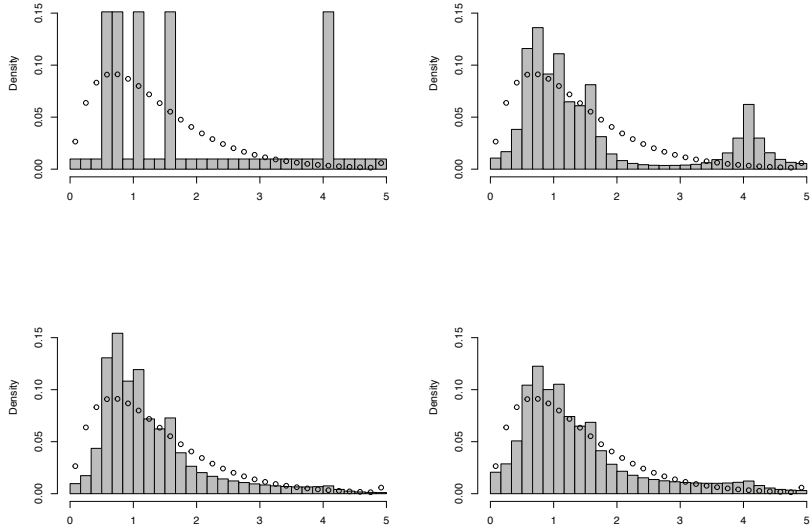}
\end{center}
\caption{Gamma$(2,1.5)$ population pmf given by dotted curve; $5$ samples across $30$ bins in range [0,8] (not pictured) with $1$ point in bins $1,2,3,7,16$.  Posterior mean mass functions from MSB$(G)$ priors:  Top left: $G_1$=Dirichlet$(w,\ldots, w)$, $w=2/29$; Bottom left: $G_2=$Tridiagonal with $w=8$; Bottom right: $G_3=$Tridiagonal with $w=16$; Top right $G_4= (G_1 + 2.5G_2)/3.5$.}
\end{figure}

\begin{figure}
\label{Fig4}
\begin{center}
\includegraphics
{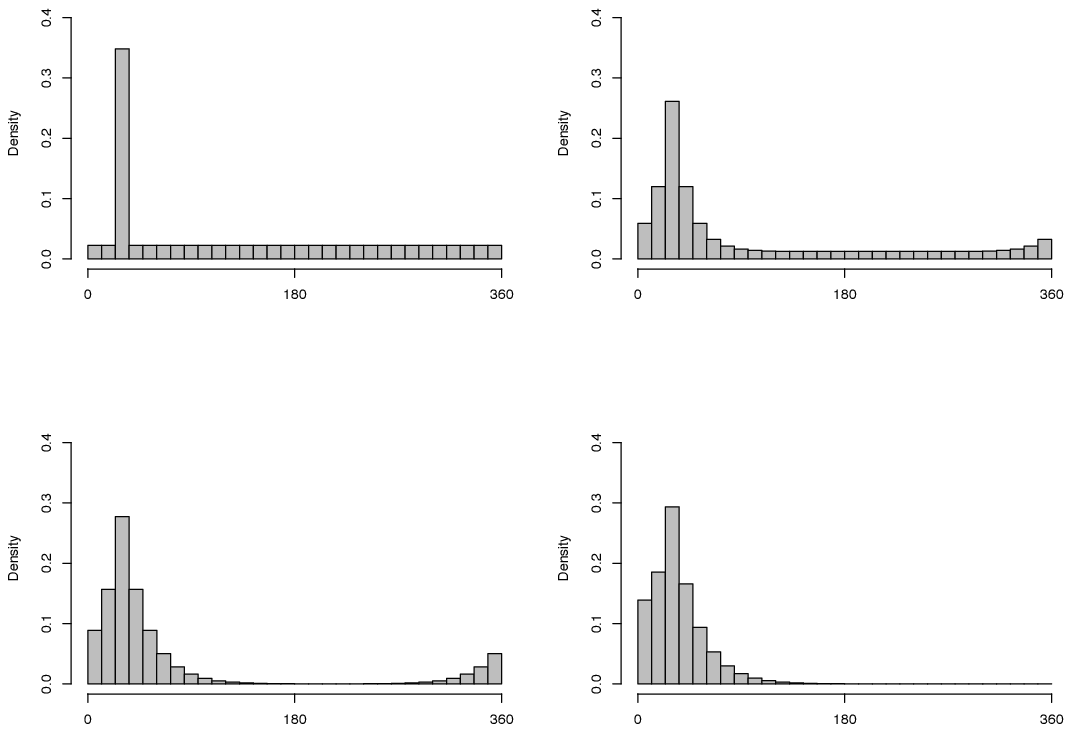}
\end{center}
\caption{Wrapped vs unwrapped; $1$ sample, across $30$ bins in degree range [0, 360] (not pictured), in bin $3$.  Posterior mean mass functions from MSG$(G)$ priors:  Top left: $G_1$=Dirichlet$(w,\ldots, w)$, $w=2/29$; Bottom left: $G_2=$Wrapped tridiagonal with $w=3$; Bottom right: $G_3=$(unwrapped) Tridiagonal with $w=3$; Top right $G_4= (G_1 + 2.5G_2)/3.5$.}
\end{figure}

\medskip \noindent
{\bf Types of generators and associations.}  We now consider several ways, among others, in which a network or graph might be specified and an associated generator matrix $G$ constructed.

In the context of this paper, graphs are connected, weighted, directed or undirected, and without self-loops. Weights should be nonnegative and the sum of weights of edges connected to (undirected) or coming into (directed) any one edge should have finite upper bound. 

In general, once a graph has been specified, a generator matrix $G$ is obtained from the adjacency matrix $A$ of the graph by modifying the diagonal entries of $A$ to give it a generator matrix structure. Note then that connectedness of the graph would ensure irreducibility of $G$.
In the case of infinitely many categories, we would further demand that a graph result in a positive recurrent generator $G$.

\vskip .1cm
	\textit{Dirichlet graphs} The Dirichlet prior is a special case of the Markovian stick-breaking prior. For the purpose of comparison, we begin by specifying the graph or network associated with a Dirichlet$(\vec \alpha)$ prior on $d$ categories. The corresponding graph on $d$ nodes has a directed edge from node $i$ to node $j$ of weight $\alpha_j$ for each ordered pair of distinct nodes $(i,j)$. Thus, for every node $j$, all incoming edges have weight $\alpha_j$ independent of the originating node, disallowing for special associations between pairs of nodes. The adjacency matrix $A$ for this graph is constant with $A_{ij}=\alpha_j$, and the associated generator matrix $G$ has the same off diagonal entries and diagonal entry $G_{jj}=\alpha_j-\sum_{i=1}^d\alpha_i$.
\vskip .1cm
	
	\textit{Geometric graphs} When categorical data arise from binning continuous data, categories come with a geometric arrangement. For ease, suppose the continuous data is real-valued data, and so categories (intervals in which the continuous data occur) come linearly ordered. This geometric arrangement can be reflected in a graph with categories represented by nodes and an undirected edge of weight $w$ placed between each pair of adjacent categories, forming a line segment. The adjacency matrix $A$ for such a graph has $w$ in the first upper diagonal and first lower diagonal entries and zeros elsewhere. The associated generator matrix $G$ has the same off diagonal entries as $A$ and the necessary diagonal entries for generator structure. We will refer to this type of generator $G$ as `tridiagonal' with weight $w$. We mention that the prior MSB$(G)$ mean, in this case, would be uniform.  Moreover, the weight $w$ represents a relative strength, and can be related to $\theta$ when $G$ is put in form $G=\theta(Q-I)$.  By increasing $w$, the `smoothing' effect, relative to the geometry, on the posterior mean estimate of the pmf will strengthen.

	There are of course other relevant settings.  For instance, suppose the continuous data were angle data taking values on the circle and having full support. In such a case, the corresponding graph would be a cycle graph on $d$ categories, and the corresponding `wrapped' generator matrix would be obtained by modifying the tridiagonal generator with weight $w$ to have entries $G_{d,1}=G_{1,d}=w$ and $G_{11}=G_{dd}=-2w$.

	More complicated geometries can also be envisioned, for instance when the categories of interest are regions in a mesh of a many-dimensional setting. 
	
	\vskip .1cm
	{\it Contingency tables.}
	The Markovian stick-breaking prior might also be used for multi-factor categorical data, where a single data point is of the form $x=(s_1,s_2,...,s_k)\in \XX=S_1\times S_2\times\cdots\times S_k$, representing $k$ categorical factors observed, where $S_i$ is the set of possible outcomes of factor $i$. As an example, one might simultaneously observe eye and hair color of individuals. Then $k=2$ and $S=\{$eye colors$\}\times\{$hair colors$\}$, and a single observation might be (brown eyes, black hair). In certain contexts, such as genetics, we might have prior reason to believe that similar outcomes (differing by only a few factors) are similarly likely to occur in the population. Thus, a prior distribution on $\Delta_\XX$ should put more weight on distributions where similar outcomes have similar probabilities of occurring.

In specifying such an MSB prior on $\Delta_\XX$, we might translate the notion of similar outcomes into a network. For example: For two outcomes $x=(s_1,s_2,...,s_k)$ and $y=(t_1,t_2,...,t_k)$, place an undirected edge of weight $w$ between them only if the outcomes are identical for all but one factor $s_j\neq t_j$. Such a prior associates any two outcomes differing only by a single factor. Interestingly, as the associated generator matrix is by construction lumpable according to each factor, this joint MSB prior on $\Delta_\XX$ has marginal Dirichlet$(w,w,...,w)$ prior on $\Delta_{S_i}$ for each factor. Similarly, we might specify a joint prior on $\Delta_{\XX}$ with pre-specified MSB$(G^{(i)})$ marginals on each $\Delta_{S_i}$ which encodes closeness of similar outcomes in $\XX$ by defining a joint generator matrix $G_{x,y}=G_{s_j,t_j}^{(j)}$ for $x=(s_1,s_2,...,s_k)$ and $y=(t_1,t_2,...,t_k)$ identical for all but one factor $s_j\neq t_j$, and $G_{x,y}=0$ otherwise for $x\neq y$.

	 \vskip .1cm
	 
	\textit{Directed vs undirected graphs.}  Since a graph with undirected edges corresponds to a symmetric adjacency matrix, the associated Markovian stick-breaking prior will correspond to a symmetric $G$ with a uniform stationary vector.  Necessarily then, an MSB, with non-uniform prior mean vector, corresponds to a directed graph.  Note that some directed graphs also produce a uniform mean stationary vector, such as a directed cycle graph with equal weights.

One might envision using directed graphs in settings where there is a `hierarchy', such as in employee data in different levels of management, for instance.

\medskip \noindent

{\bf Simple numerical experiments.}
In Theorem \ref{thmpriornset}, we have computed the posterior mean estimate of the probability mass function given that the prior is a Markovian stick-breaking measure with generator $G$ and the empirical counts $\vec{k}$ of observed data.  Specifically, let $\XX$ denote the set of categories and let $(\nu,T_1)\sim$ MSB$(G)$, where $\nu$ is the prior. For a data vector $k=(k_w)_{w\in\XX}$ of non-negative integers and a category $x\in\XX$, define $v(k)=\big(v_x(k)\big)_{x\in\XX}$ by
	$$v_x(k)=\E\left[\prod_{w\in\XX}\nu(w)^{k_w}\Big|T_1=x\right]=\binom nk^{-1}\sum_{\sigma\in\mathbb S(k)}\left(I-\frac Gn\right)^{-1}_{x,\sigma_n}\prod_{j=1}^{n-1}\left(I-\frac Gj\right)^{-1}_{\sigma_{j+1},\sigma_j}$$
	
	where $n=\sum_{w\in\XX}k_w$ and $\mathbb S(k)$ denotes the set of distinct permutations of a list containing precisely $k_w$ many $w$'s for each $w\in\XX$. Then, the posterior probability mass function, specified in terms of the posterior means, given the observed multinomial counts $k$, when evaluated at $x\in\XX$, is given by
	$$p(x|k)=\frac{\mu^Tv(k+e_x)}{\mu^Tv(k)}=\frac{\sum_{w\in\XX}\mu_wv_w(k+e_x)}{\sum_{w\in\XX}\mu_wv_w(k)}$$
	where $\mu$ is the stationary vector of $G$.

We consider now simple computational experiments to see how different generators $G$, with respect to Markovian stick-breaking priors, affect the posterior mean probability mass function, computed exactly from the above formulas with a small number of samples, in two types of data, one with Normal and the other with Gamma samples.  We will consider $G$'s, which are Dirichlet, tri-diagonal, and averages between these types, to see the effects.

In Figure 2,
with respect to a generated Normal$(0,1)$ sample histogram of $6$ samples, across $30$ bins from $-5$ to $5$, with 
$1$ point in bins $10$ and $12$ and $2$ points in bins $15$ and $17$, posterior mean mass functions are plotted with respect to four Markovian stick-breaking priors.  
In the top left plot, the generator $G_1$ corresponds to a Dirichlet$(w,\ldots, w)$ where $w=2/29$.  In the bottom left and bottom right, the generators $G_2$ and $G_3$ are a tridiagonal matrices with $w=3$ and $w=8$ entries in the two off-diagonals respectively.  In the top right, the generator $G_4$ is the average $G_4=(G_1+2.5G_2)/(3.5)$. 

Similarly, in Figure 3,
with respect to a generated Gamma$(2,1.5)$ sample histogram of $5$ samples, again across $30$ bins from $0$ to $8$, with $1$ point in bins $1,2,3,7$ and $16$, posterior mean mass functions are plotted with respect to similar priors in the same locations as in Figure 2.

In Figure 4, the intent is to see the posterior mean mass function effects, with respect to one data point in bin $3$, across $30$ bins indexed by angles (degrees) of a circle, when the priors correspond to generators which are wrapped tri-diagonal $G_2$ with $w=3$ in the bottom left, $G_1=$Dirichlet$(w, \ldots, w)$ with $w=2/29$ in the top left, their average $G_4= (G_1+2.5G_2)/(3.5)$ in the top right, and an unwrapped tri-diagonal generator $G_3$ with $w=3$ in the bottom right. 

\medskip
{\it Discussion.}
Briefly, we were interested to see what effects might arise from using Markovian stick-breaking priors in probability mass function inference.   We observe in Figures 2 and 3 that the posterior mean mass functions, computed from Markovian stick-breaking priors with tridiagonal $G$'s in the bottom left and right, show clear effects due to the network affinities encoded in the generators in comparison to the posterior mean mass function with respect to the
 Dirichlet prior in the top left.  The posterior mean mass function with respect to the prior built with the averaged $G_4$ generator incorporates a similarity structure with some positive weight between all categories, but with emphasis on neighbor categories.  In Figure 4, one definitely sees the effect of wrapping in the bottom left, and also the averaging effect where all bins receive non-negligible mass in the top right.
 
  It would seem that similarities between categories encoded in the generator $G$ do affect the posterior distribution when the prior is a Markovian stick-breaking process with generator $G$.  
 In terms of future work, there are of course several natural directions to pursue, among them to clarify more the scope and performance of these Markovian stick-breaking priors in various categorical network settings.

\section{Proof of Theorem \ref{thmpriornset}}
\label{proof_sec}

We begin by enumerating some facts.

\vskip .1cm
 {\it Fact 1.}
	 Let $\P\sim$ GEM$(\theta)$. Then $P_1\sim$ Beta$(1,\theta)$ and
	\begin{align}
	\E\big[(1-P_1)^jP_1^{k-j}\big] & 
	=\frac{\Gamma(1+\theta)}{\Gamma(1)\Gamma(\theta)}\frac{\Gamma(1+k-j)\Gamma(\theta+j)}{\Gamma(1+\theta+k)}
	=\frac{\theta\Gamma(k-j+1)\Gamma(\theta+j)}{\Gamma(\theta+k+1)}
	\nonumber\\
	\E\big[(1-P_1)^k\big] & 
	=\frac{\theta\Gamma(1)\Gamma(\theta+k)}{\Gamma(\theta+k+1)}
	=\frac\theta{\theta+k}
	\label{factmoment}
	\end{align}
	
	\vskip .1cm
	
	{\it Fact 2.} 
	Let $G=\theta(Q-I)$ be a generator matrix and $Q$ stochastic and $\theta>0$. When $k>0$:
	\begin{equation}
	\left(I-\frac{\theta Q}{\theta+k}\right)^{-1}=\frac{\theta+k}k(I-G/k)^{-1} \ {\rm and \ } Q(I-G/k)^{-1}=\frac{\theta+k}\theta\left((I-G/k)^{-1}-\frac{kI}{\theta+k}\right)
	\label{factqalg} \end{equation}
	
	\vskip .1cm
	{\it Fact 3.} Consider the space $\{0,1,2,\ldots\}^n$ of non-negative integer $n$-vectors. For two vectors $\vec k,\vec l\in\{0,1,2,\ldots\}^n$, we say $\vec l<\vec k$ if for each $1\le j\le n$, we have $l_j\le k_j$, and for some $1\le j\le n$, in fact $l_j<k_j$. Note that this gives a strict partial ordering to all non-negative $n$-vectors; that the zero vector is strictly less than every other vector; and that each $\vec k$ is strictly greater than only finitely many $n$-vectors. Thus, for each $n$, the space is well-founded and an induction may be considered with respect to this partial ordering starting from 0.

	\vskip .1cm
	
	 {\it Fact 4.} For an $n$-vector $\vec k$ of non-negative integer entries, with $k=\sum_{j=1}^nk_j>0$,
	\begin{equation}
	\#\mathbb{S}(\vec k)=\binom k{k_1,k_2,...,k_n}=\frac{\Gamma(k+1)}{\prod_{j=1}^n\Gamma(k_i+1)}
	\label{factS} \end{equation}

The following proposition will help an induction in the proof of Theorem \ref{thmpriornset}.
\begin{prop}
	Let $G$ be an irreducible, positive recurrent generator matrix on $\XX$ and let $(\nu,T_1)\sim$ MSB$(G)$. Then, for each $k\in\{0,1,2,\ldots\}$, $A\subseteq\XX$, and $x\in\XX$, we have
	\begin{align}
	\label{prop_help}
	\E\Big[\nu(A)^k\Big|T_1=x\Big]=e_x^T\prod_{j=1}^k\Big((I-G/j)^{-1}D(A)\Big)\vec 1
	\end{align}
	\label{propprioroneset}
\end{prop}
\vspace{-4mm}
\begin{proof}
	Since \eqref{prop_help} is a statement regarding the distribution of $(\nu,T_1)$, we can choose a particular instance of $(\nu,T_1)$ constructed from an independent pair $\X=(X_j)_{j=1}^\infty$ and $\T=(T_j)_{j=1}^\infty$ of, respectively, an i.i.d. sequence of Beta$(1,\theta)$ variables and a stationary, homogeneous Markov chain with transition kernel $Q$, where $G=\theta(Q-I)$. 
	
	As usual, let $\P$ be defined with respect to $\X$ by $P_j=X_j\prod_{i=1}^{j-1}(1-X_i)$.
	For ease of notation, define the vector $v(k,A)=\big(v_x(k,A)\big)_{x\in\XX}$ by
	$v_x(k,A)=\E\big[\nu(A)^k\big|T_1=x\big]$.
	We begin by finding a recursive (in $k$) formula for $v(k,A)$.

	 To this end, we define
	$\nu^*=\sum_{j=2}^\infty \left[X_j\prod_{i=2}^{j-1}(1-X_i)\right]\delta_{T_j}$
	and note that $\nu^*$ is independent of $X_1=P_1$ since $\X$ is i.i.d. and independent of $\T$. Furthermore,
	$\nu=P_1\delta_{T_1}-(1-P_1)\nu^*$.
Write	
	\begin{align}
		& e_x^Tv(k,A)
		 \ = \ v_x(k,A) \ = \ 
		\E\Big[\left(P_1\delta_x(A)+(1-P_1)\nu^*(A)\right)^k\big|T_1=x\Big]\nonumber
		\\
		& =\sum_{y\in\XX}\p\big(T_2=y|T_1=x\big)\E\Big[\left(P_1\delta_x(A)+(1-P_1)\nu^*(A)\right)^k\big|T_1=x,T_2=y\Big]\nonumber
		\\
		& =\sum_{y\in\XX}Q_{xy}\Bigg\{\E\Big[\big((1-P_1)\nu^*(A)\big)^k\big|T_1=x,T_2=y\Big]\nonumber \\
		&\ \ \ \ \ \ \ \ \  +\delta_x(A)\sum_{j=0}^{k-1}\binom{k}{j}\E\Big[P_1^{k-j}(1-P_1)^j\big(\nu^*(A)\big)^j\big|T_1=x,T_2=y\Big]\Bigg\}
		\label{eqa}
	\end{align}
	Clearly, as $\X$ and $\T$ are independent and $\X$ is i.i.d., $P_1=X_1$ is independent of $T_1$, $T_2$, and $\nu^*$. By the Markov property and since $\nu^*$ is not a function of $T_1$, we have $\nu^*|(T_2=y)$ is independent of $T_1$. Furthermore, since $(X_j)_{j\ge1}\stackrel d=(X_j)_{j\ge2}$ as an i.i.d. sequence and $(T_j)_{j\ge1}\stackrel d=(T_j)_{j\ge2}$ as a stationary Markov chain, we have $(\nu,T_1)\stackrel d=(\nu^*,T_2)$, implying $\nu^*|(T_1=x,T_2=y)\stackrel d=\nu|(T_1=y)$. 
	
	Thus, defining $\nu(0,A)=\vec 1$, equation \eqref{eqa} becomes
	\begin{align*}
		 & =\sum_{y\in\XX}Q_{xy}\hspace{-1.2mm}\left[\E\left[(1-P_1)^k\right]\E\left[\nu(A)^k\big|T_1=y\right]\hspace{-.4mm}+\hspace{-.2mm}\delta_x(A)\sum_{j=0}^{k-1}\binom{k}{j}\E\big[P_1^{k-j}(1-P_1)^j\big]\E\hspace{-.4mm}\left[\nu(A)^j\big|T_1=y\right]\hspace{-.9mm}\right]\nonumber
		 \\
		 & =\sum_{y\in\XX}Q_{xy}\left[\E\left[(1-P_1)^k\right]v_y(k,A)+\delta_x(A)\sum_{j=0}^{k-1}\binom{k}{j}\E\left[P_1^{k-j}(1-P_1)^j\right]v_y(j,A)\right]\\
		 &= e_x^T\left[\E\left[(1-P_1)^k\right]Qv(k,A)+D(A)Q\sum_{j=0}^{k-1}\binom{k}{j}\E\left[P_1^{k-j}(1-P_1)^j\right]v(j,A)\right],
		 \end{align*}
		 which, noting \eqref{factmoment}, equals
		 \begin{align*}
		 &=e_x^T\left[\frac\theta{\theta+k}Qv(k,A)+D(A)Q\sum_{j=0}^{k-1}\binom{k}{j}\frac{\theta\Gamma(k-j+1)\Gamma(\theta+j)}{\Gamma(\theta+k+1)}v(j,A)\right].
	\end{align*}

	Since the statement holds for every $x$, it follows that 
	\begin{align*}
		v(k,A) & =\frac\theta{\theta+k}Qv(k,A)+D(A)Q\sum_{j=0}^{k-1}\binom{k}{j}\frac{\theta\Gamma(k-j+1)\Gamma(\theta+j)}{\Gamma(\theta+k+1)}v(j,A)\nonumber
		\\
		& =\frac\theta{\theta+k}Qv(k,A)+\frac{\theta\Gamma(k+1)}{\Gamma(\theta+k+1)}D(A)Q\sum_{j=0}^{k-1}\frac{\Gamma(\theta+j)}{\Gamma(j+1)}v(j,A)\end{align*}
and		
		$\left(I-\frac{\theta Q}{\theta+k}\right)v(k,A)  =\frac{\theta\Gamma(k+1)}{\Gamma(\theta+k+1)}D(A)Q\sum_{j=0}^{k-1}\frac{\Gamma(\theta+j)}{\Gamma(j+1)}v(j,A)$.
		Then,
		\begin{align}
		v(k,A) & =\frac{\theta\Gamma(k+1)}{\Gamma(\theta+k+1)}\left(I-\frac{\theta Q}{\theta+k}\right)^{-1}D(A)Q\sum_{j=0}^{k-1}\frac{\Gamma(\theta+j)}{\Gamma(j+1)}v(j,A)\nonumber
		\\
		& =\frac{\theta\Gamma(k)}{\Gamma(\theta+k)}\left(I-\frac Gk\right)^{-1}D(A)Q\sum_{j=0}^{k-1}\frac{\Gamma(\theta+j)}{\Gamma(j+1)}v(j,A),\label{eqpriorrecursion}
	\end{align}
	where the last line follows from \eqref{factqalg}.
	
	We now solve the recursion
	for $v(k,A)$ inductively. We have already specified $v(0,A)=\vec 1$. By \eqref{eqpriorrecursion}, we have
	$$v(1,A)=\frac{\theta\Gamma(1)}{\Gamma(\theta+1)}\left(I-G/1\right)^{-1}D(A)Q\frac{\Gamma(\theta)}{\Gamma(1)}v(0,A)=\left(I-G/1\right)^{-1}D(A)Q\vec 1=\left(I-G/1\right)^{-1}D(A)\vec 1$$
	If, for $1\le j\le k-1$, 
	$v(j,A)=\prod_{i=1}^j\Big((I-G/i)^{-1}D(A)\Big)\vec 1$,
	then it follows from \eqref{eqpriorrecursion} that
	$v(k,A) = \frac{\theta \Gamma(k)}{\Gamma(\theta+k)}\big(I-G/k\big)^{-1} u_k$
	where
	$$u_k = D(A)Q\sum_{j=0}^{k-1}\frac{\Gamma(\theta+j)}{\Gamma(j+1)}\prod_{i=1}^j\Big((I-G/i)^{-1}D(A)\Big)\vec 1.$$
	
	We now claim that $u_k=w_k$ where
	\begin{align*}
	w_k = \frac{\Gamma(\theta+k)}{\theta\Gamma(k)}D(A)\prod_{i=1}^{k-1}\Big((I-G/i)^{-1}D(A)\Big)\vec 1.
	\end{align*}
	Indeed, if $u_k = w_k$, we would conclude that
	$$v(k,A) = \frac{\theta\Gamma(k)}{\Gamma(\theta+k)}\left(I-\frac Gk\right)^{-1}w_k = \prod_{i=1}^k\Big((I-G/i)^{-1}D(A)\Big)\vec 1,$$
	finishing the proof of Proposition \ref{propprioroneset}.
	
	To verify the claim, observe that $u_1 = \Gamma(\theta)D(A)Q\vec 1 = \Gamma(\theta)D(A)\vec 1 = \big(\Gamma(\theta +1)/\theta\big)D(A)\vec 1 = w_1$.  Suppose that $u_j = w_j$ for $j\leq k$.  Then,
	\begin{align*}
	u_{k+1} - w_{k+1} &= u_k + \frac{\theta}{k}D(A)Q(I-G/k)^{-1}w_k - \frac{\theta+k}{k}D(A)(I-G/k)^{-1}w_k\\
	&= u_k -D(A)w_k = u_k-w_k = 0,
	\end{align*}
as $(D(A))^2=D(A)$, finishing the proof.  \end{proof}	

 \noindent
{\bf Proof of Theorem \ref{thmpriornset}.}	
	The theorem holds trivially for $\vec k=\vec 0$. If $\vec k\neq\vec 0$, without loss of generalization, we assume $\vec k$ has strictly positive entries. Otherwise, it may be represented as a vector $\vec k^\sharp$ of smaller length by omitting the zero entries, with $\vec A^\sharp$ the corresponding shortened vector of sets and $n^\sharp$ the new vector length.
	
	As in the proof of the Proposition \ref{propprioroneset}, 
	 we may choose a particular instance of $(\nu,T_1)$ constructed from an independent pair $\X=(X_j)_{j=1}^\infty$ and $\T=(T_j)_{j=1}^\infty$ of, respectively, an i.i.d. sequence of Beta$(1,\theta)$ variables and a stationary, homogeneous Markov chain with transition kernel $Q$, where $G=\theta(Q-I)$. 	Let $\P$ be defined with respect to $\X$ by $P_j=X_j\prod_{i=1}^{j-1}(1-X_i)$.
	
	Define now the vector $v(\vec k,\vec A)=\big(v_x(\vec k,\vec A)\big)_{x\in\XX}$ by
	$v_x(\vec k,\vec A)=\E\big[\prod_{j=1}^n\nu(A_j)^{k_j}\big|T_1=x\big]$.
	We begin by finding a recursive (in $\vec k$ and $n$) formula for $v(\vec k,\vec A)$, and then we solve the recursion using Lemma \ref{lemmapriorinduction2}, stated at the end of the section. 
	
	To this end, recall the definition
	$\nu^*=\sum_{j=2}^\infty \big[X_j\prod_{i=2}^{j-1}(1-X_i)\big]\delta_{T_j}$.
	We compute
	\begin{align*}
		& e_x^Tv(\vec k,\vec A)
		\  = \ \E\left[\prod_{j=1}^n\left(P_1\delta_{T_1}(A_j)+(1-P_1)\nu^*(A_j)\right)^{k_j}\Big|T_1=x\right] \nonumber
		\\
		& =\sum_{y\in\XX}\p\big(T_2=y|T_1=x\big) \E\left[\prod_{j=1}^n\left(P_1\delta_x(A_j)+(1-P_1)\nu^*(A_j)\right)^{k_j}\Big| T_1=x,T_2=y\right], \nonumber 
	\end{align*}
		which equals, as the collection $\vec A$ consists of disjoint set so that $\delta_x(A_i)\delta_x(A_j) =0$ for $i\neq j$,
\begin{align}
\label{star}		
		&\sum_{y\in\XX}Q_{x,y}\left\{\E\left[\prod_{j=1}^n\left((1-P_1)\nu^*(A_j)\right)^{k_j}\Big| T_1=x,T_2=y\right]\right. \\
		& \hspace{4mm}+\left.\sum_{i=1}^n\delta_x(A_i)\sum_{l=0}^{k_i-1}\binom{k_i}l\E\left[(1-P_1)^lP_1^{k_i-l}\nu^*(A_i)^l\hspace{-.9mm}\prod_{1\le j\le n;\ j\neq i}\hspace{-.9mm}\left((1-P_1)\nu^*(A_j)\right)^{k_j}\Big| T_1=x,T_2=y\right]\hspace{-.5mm}\right\}.\nonumber
	\end{align}
	 
	Recall the relations among $P_1$, $\nu$, $\nu^*$, $T_1$ and $T_2$ stated below \eqref{eqa}. Then, by Fact 1, we have that \eqref{star} equals
	\begin{align}
		&\sum_{y\in\XX}Q_{x,y}\left[\E\left[(1-P_1)^k\right]\E\left[\prod_{j=1}^n\nu(A_j)^{k_j}\Big|T_1=y\right]+\right.\nonumber\\
		& \hspace{5mm}+\left.\sum_{i=1}^n\delta_x(A_i)\sum_{l=0}^{k_i-1}\binom{k_i}l\E\left[(1-P_1)^{k-k_i+l}P_1^{k_i-l}\right]\E\left[\nu(A_i)^l\prod_{1\le j\le n;\ j\neq i}\left(\nu(A_j)\right)^{k_j}\Big| T_1=y\right]\right]\nonumber
		\\
		& =\sum_{y\in\XX}Q_{x,y}\Bigg[\frac\theta{\theta+k}v_y(\vec k,\vec A)+\nonumber\\
		& \hspace{5mm}+\left.\sum_{i=1}^n\delta_x(A_i)\sum_{l=0}^{k_i-1}\binom{k_i}l\frac{\theta\Gamma(k_i-l+1)\Gamma(\theta+k-k_i+l)}{\Gamma(\theta+k+1)}v_y(\vec k+(l-k_i)e_i,\vec A)\right]\nonumber
		\\
		& =e_x^T\left[\frac\theta{\theta+k}Qv(\vec k,\vec A)+\sum_{i=1}^n\frac{\theta\Gamma(k_i+1)}{\Gamma(\theta+k+1)}D(A_i)Q\sum_{l=0}^{k_i-1}\frac{\Gamma(\theta+k-k_i+l)}{\Gamma(l+1)}v(\vec k+(l-k_i)e_i,\vec A)\right]\nonumber
	\end{align}
	
	Since the above computation holds for every $x$, it may be written as a vector equation:
	\begin{align}
		v(\vec k,\vec A) & =\frac\theta{\theta+k}Qv(\vec k,\vec A)+\sum_{i=1}^n\frac{\theta\Gamma(k_i+1)}{\Gamma(\theta+k+1)}D(A_i)Q\sum_{l=0}^{k_i-1}\frac{\Gamma(\theta+k-k_i+l)}{\Gamma(l+1)}v(\vec k+(l-k_j)e_j,\vec A)\nonumber\\
		& =\left(I-\frac{\theta Q}{\theta+k}\right)^{-1}\sum_{i=1}^n\frac{\theta\Gamma(k_i+1)}{\Gamma(\theta+k+1)}D(A_i)Q\sum_{l=0}^{k_i-1}\frac{\Gamma(\theta+k-k_i+l)}{\Gamma(l+1)}v(\vec k+(l-k_j)e_j,\vec A)\nonumber
		\\
		& =\left(I-G/k\right)^{-1}\sum_{i=1}^n\frac{\theta\Gamma(k_i+1)}{k\Gamma(\theta+k)}D(A_i)Q\sum_{l=0}^{k_i-1}\frac{\Gamma(\theta+k-k_i+l)}{\Gamma(l+1)}v(\vec k+(l-k_i)e_i,\vec A).\label{eqvrecursive}
	\end{align}

	The recursive formula \eqref{eqvrecursive} for $v(\vec k,\vec A)$ is in terms of the values of $v(\vec l,\vec A)$ only for $\vec l<\vec k$. If $\vec k$ has $r$ zero entries and at least one positive entry, recall $\vec k^\sharp$, the reduction of $\vec k$ to a strictly positive $(n-r)$-vector by removal of zero entries, with $\vec A^\sharp$ corresponding. Then $v(\vec k,\vec A)=v(\vec k^\sharp,\vec A^\sharp)$.
	Thus, we consider simultaneously an induction on the value of $n$ and, given $n$, an induction on the $n$-vector $\vec k$ according to the strict partial ordering from Fact 3.

	When $n=1$, the theorem holds by Proposition \ref{propprioroneset}. This is the base case for induction on $n$. Suppose by way of induction on $n$ that, for each $1\le m<n$ and, given $m$, each non-negative integer $m$-vector $\vec l$ and $m$-vector $\vec B$ of disjoint subsets of $\XX$, we have that the theorem holds for $v(\vec l,\vec B)$.
	
	Consider a non-negative integer $n$-vector $\vec k$ with at least one positive entry and an $n$-vector $\vec A$ of disjoint subsets of $\XX$. If $\vec k$ has any zero-entries, by the induction assumption on $n$, 
	$v(\vec k,\vec A)=v(\vec k^\sharp,\vec A^\sharp)=\left(\#\mathbb{S}(\vec k)\right)^{-1}\sum_{\sigma\in\mathbb{S}(\vec k)}\left[\prod_{r=1}^k\left((I-G/r)^{-1}D(A_{\sigma(j)})\right)\right]\vec1$.
	This is the base case for an induction on $\vec k\in\{0,1,2,\ldots\}^n$.
	
	 Suppose instead that $\vec k$ consists of positive integers. Given $n$, suppose by way of induction on $\vec k\in\{0,1,2,\ldots\}^n$ that for every $\vec l\in\{0,1,2,\ldots\}^n$ with $\vec l<\vec k$, the theorem holds for $v(\vec l,\vec A)$.
	Then, we have
	\begin{align*}
		v(\vec k,\vec A) & =\frac\theta{k\Gamma(\theta+k)}\left(I-G/k\right)^{-1}\sum_{i=1}^n\Gamma(k_i+1)D(A_i)Q\sum_{l=0}^{k_i-1}\frac{\Gamma(\theta+k-k_i+l)}{\Gamma(l+1)}v(\vec k+(l-k_i)e_i,\vec A)
		\end{align*}
		equals, using  \eqref{eqvrecursive},
		\begin{align}
		\label{star2}
		& \frac\theta{k\Gamma(\theta+k)}\left(I-G/k\right)^{-1}\sum_{i=1}^n\Gamma(k_i+1)D(A_i)Q\sum_{l=0}^{k_i-1}\frac{\Gamma(\theta+k-k_i+l)}{\Gamma(l+1)}\\
		& \hspace{15mm}\times \left(\#\mathbb{S}(\vec k+(l-k_i)e_i)\right)^{-1}\sum_{\sigma\in\mathbb{S}(\vec k+(l-k_i)e_i)}\left[\prod_{r=1}^{k-k_i+l}\left((I-G/r)^{-1}D(A_{\sigma(j)})\right)\right]\vec1\nonumber
	\end{align}
	Recalling \eqref{factS}, it then follows that \eqref{star2} equals
	\begin{align}
		& \frac\theta{k\Gamma(\theta+k)}\left(I-G/k\right)^{-1}\sum_{i=1}^n\Gamma(k_i+1)D(A_i)Q\sum_{l=0}^{k_i-1}\frac{\Gamma(\theta+k-k_i+l)}{\Gamma(l+1)}\nonumber\\
		& \hspace{15mm}\times\frac{\Gamma(l+1)\prod_{j\neq i}\Gamma(k_j+1)}{\Gamma(k-k_i+l+1)}\sum_{\sigma\in\mathbb{S}(\vec k+(l-k_i)e_i)}\left[\prod_{r=1}^{k-k_i+l}\left((I-G/r)^{-1}D(A_{\sigma(j)})\right)\right]\vec1\nonumber
		\\
		& =\frac{\prod_{j=1}^n\Gamma(k_j+1)}{k\Gamma(\theta+k)}\left(I-G/k\right)^{-1}\sum_{i=1}^n\theta D(A_i)Q\sum_{l=0}^{k_i-1}\frac{\Gamma(\theta+k-k_i+l)}{\Gamma(k-k_i+l+1)}\nonumber\\
		& \hspace{15mm}\times \sum_{\sigma\in\mathbb{S}(\vec k+(l-k_i)e_i)}\left[\prod_{r=1}^{k-k_i+l}\left((I-G/r)^{-1}D(A_{\sigma(j)})\right)\right]\vec1.
		\label{eqbigind}\end{align}
		By Lemma \ref{lemmapriorinduction2}, at the end of the section, \eqref{eqbigind} equals
		\begin{align*}
		& \frac{\prod_{j=1}^n\Gamma(k_j+1)}{k\Gamma(\theta+k)}\left(I-G/k\right)^{-1}\sum_{i=1}^n\frac{\Gamma(\theta+k)}{\Gamma(k)}D(A_i)\sum_{\sigma\in\mathbb{S}(\vec k-e_i)}\left[\prod_{j=1}^{k-1}\left((I-G/r)^{-1}D(A_{\sigma(j)})\right)\right]\vec 1\nonumber
		\\
		& =\frac{\prod_{j=1}^n\Gamma(k_j+1)}{\Gamma(k+1)}\sum_{\sigma\in\mathbb{S}(\vec k)}\left[\prod_{j=1}^k\left((I-G/r)^{-1}D(A_{\sigma(j)})\right)\right]\vec 1\nonumber
		\\
		& =\left(\#\mathbb{S}(\vec k)\right)^{-1}\sum_{\sigma\in\mathbb{S}(\vec k)}\left[\prod_{j=1}^k\left((I-G/r)^{-1}D(A_{\sigma(j)})\right)\right]\vec 1.\nonumber
	\end{align*}
	By induction on $\vec k$, the statement of the theorem holds for all $\vec k\in\{0,1,2,\ldots\}^n$. By induction on $n$, the theorem holds for all $n$ as well. This completes the proof. 
\qed

\medskip
We now state and prove the lemma referred to in the argument for 
Theorem \ref{thmpriornset}.

\begin{lemma}
	Let $m\ge 1$, $(A_j)_1^n$ be disjoint sets, $\vec k\in \{0,1,2,\ldots\}^{n-1}$, and $\tilde k=\sum_{j=1}^{n-1}k_j$. Then,
	\begin{align}
	\label{star3}
		 \theta D(A_n)Q&\sum_{l=0}^{m-1}\frac{\Gamma(\theta+\tilde k+l)}{\Gamma(\tilde k+l+1)}\sum_{\sigma\in\mathbb{S}(\vec k,l)}\prod_{j=1}^{\tilde k+l}\left[(I-G/j)^{-1}D(A_{\sigma(j)})\right]\vec 1\\
		= & \frac{\Gamma(\theta+\tilde k+m)}{\Gamma(\tilde k+m)}D(A_n)\sum_{\sigma\in\mathbb{S}(\vec k,m-1)}\prod_{j=1}^{\tilde k+m-1}\left[(I-G/j)^{-1}D(A_{\sigma(j)})\right]\vec 1\nonumber
	\end{align}
	where $\mathbb S(\vec k,l)=\mathbb S((k_1,...,k_{n-1},l))$.
	\label{lemmapriorinduction2}
\end{lemma}

\begin{proof} We prove the lemma by induction. Define the left-hand side of \eqref{star3} as $u_m$. Then, by \eqref{factqalg} and that $D(A_i)D(A_j)=D(A_i)\delta_i(j)$, we have 
	\begin{align*}
	u_1 & =\theta D(A_n)Q\frac{\Gamma(\theta+\tilde k)}{\Gamma(\tilde k+1)}\sum_{\sigma\in\mathbb{S}(\vec k,0)}\prod_{j=1}^{\tilde k+0}\left[(I-G/j)^{-1}D(A_{\sigma(j)})\right]\vec 1\nonumber\\
	& =\frac{\theta+k}\theta\theta\frac{\Gamma(\theta+\tilde k)}{\Gamma(\tilde k+1)}D(A_n)\sum_{\sigma\in\mathbb{S}(\vec k,0)}\prod_{j=1}^{\tilde k+0}\left[(I-G/j)^{-1}D(A_{\sigma(j)})\right]\vec 1\\
	& =\frac{\Gamma(\theta+\tilde k+1)}{\Gamma(\tilde k+1)}D(A_n)\sum_{\sigma\in\mathbb{S}(\vec k,0)}\prod_{j=1}^{\tilde k+0}\left[(I-G/j)^{-1}D(A_{\sigma(j)})\right]\vec 1.\nonumber
	\end{align*}
	
	By \eqref{factqalg} and $D(A_i)D(A_j)=D(A_i)\delta_i(j)$ again, 
	\begin{align*}
	&u_{m+1}  =u_m+\theta D(A_n)Q\frac{\Gamma(\theta+\tilde k+m)}{\Gamma(\tilde k+m+1)}\sum_{\sigma\in\mathbb{S}(\vec k,m)}\prod_{j=1}^{\tilde k+m}\left[(I-G/j)^{-1}D(A_{\sigma(j)})\right]\vec 1\\
	& =u_m+\frac{\theta+\tilde k+m}\theta\theta\frac{\Gamma(\theta+\tilde k+m)}{\Gamma(\tilde k+m+1)}D(A_n)\\
	&\times
	\left[\sum_{\sigma\in\mathbb{S}(\vec k,m)}\prod_{j=1}^{\tilde k+m}\left[(I-G/j)^{-1}D(A_{\sigma(j)})\right]-\frac{\tilde k+m}{\theta+\tilde k+m}\sum_{\sigma\in\mathbb{S}(\vec k,m-1)}\prod_{j=1}^{\tilde k+m-1}\left[(I-G/j)^{-1}D(A_{\sigma(j)})\right]\right]\vec 1,
	\end{align*}
	which further equals
	\begin{align*}
	& \frac{\Gamma(\theta+\tilde k+m+1)}{\Gamma(\tilde k+m+1)}D(A_n)\sum_{\sigma\in\mathbb{S}(\vec k,m)}\prod_{j=1}^{\tilde k+m}\left[(I-G/j)^{-1}D(A_{\sigma(j)})\right]\vec 1\\
	& \hspace{15mm}+u_m-\frac{\Gamma(\theta+\tilde k+m)}{\Gamma(\tilde k+m)}D(A_n)\sum_{\sigma\in\mathbb{S}(\vec k,m-1)}\prod_{j=1}^{\tilde k+m-1}\left[(I-G/j)^{-1}D(A_{\sigma(j)})\right]\vec 1.
	\end{align*}
	We conclude the result via induction.
\end{proof}

\medskip \noindent
{\bf Acknowledgements.}  We thank J. Sethuraman for reading and comments on a draft of this manuscript.
This research was partly supported by ARO-W911NF-18-1-0311 and a Simons Foundations Sabbatical grant.





\end{document}